\documentclass[12pt]{article}
\usepackage{amssymb,latexsym}
\input epsf

\title{\textbf{Edgewise subdivisions, local $h$-polynomials and excedances in the wreath product $\ZZ_r \wr \mathfrak{S}_n$}}

\author{Christos~A.~Athanasiadis\\
Department of Mathematics\\
University of Athens\\
Athens 15784, Hellas (Greece)\\
{\small Email: \texttt{caath@math.uoa.gr}}
}

\date{{\small October 2, 2013; revised on May 12, 2014}}

  \def\NN{{\mathbb N}}

  \def\ZZ{{\mathbb Z}}
  \def\RR{{\mathbb R}}
  
  \def\raa{{\mathrm {\mathbf a}}}
  \def\ree{{\mathrm {\mathbf e}}}
  \def\ruu{{\mathrm {\mathbf u}}}
  \def\rvv{{\mathrm {\mathbf v}}}
  \def\rww{{\mathrm {\mathbf w}}}
  \def\rzz{{\mathrm {\mathbf z}}}

  \def\dD{{\mathcal D}}
  \def\eE{{\mathcal E}}
  
  \def\pP{{\mathcal P}}
  
  \def\csum{{\rm csum}}
  \def\Asc{{\rm Asc}}
  
  \def\des{{\rm des}}
  \def\ex{{\rm exc}}
  \def\fdes{{\rm fdes}}
  \def\fex{{\rm fexc}}

  \def\sd{{\rm sd}}
  \def\stab{{\rm stab}}
  \def\supp{{\rm supp}}
  \def\sm{\smallsetminus}
  \newcommand{\qed}{$\hfill \Box$}

\begin{document}
\maketitle

\newtheorem{theorem}{Theorem}[section]
\newtheorem{proposition}[theorem]{Proposition}
\newtheorem{corollary}[theorem]{Corollary}
\newtheorem{defn}[theorem]{Definition}
\newtheorem{remark}[theorem]{Remark}
\newtheorem{lemma}[theorem]{Lemma}
\newtheorem{example}[theorem]{Example}
\newtheorem{examples}[theorem]{Examples}
\newtheorem{conjecture}[theorem]{Conjecture}
\newtheorem{fact}[theorem]{Fact}
\newtheorem{question}[theorem]{Question}
\newtheorem{observation}[theorem]{Observation}
\newtheorem{claim}[theorem]{Claim}

\begin{abstract}
The coefficients of the local $h$-polynomial of the barycentric subdivision
of the simplex with $n$ vertices are known to count derangements in the
symmetric group $\mathfrak{S}_n$ by the number of excedances. A generalization of this interpretation is given for the local 
$h$-polynomial of the $r$th
edgewise subdivision of the barycentric subdivision of the simplex. This polynomial is shown to be $\gamma$-nonnegative and a combinatorial
interpretation to the corresponding $\gamma$-coefficients is provided. The new combinatorial interpretations involve the notions of flag excedance and
descent in the wreath product $\ZZ_r \wr \mathfrak{S}_n$. A related result
on the derangement polynomial for $\ZZ_r \wr \mathfrak{S}_n$, studied by Chow
and Mansour, is also derived from results of Linusson, Shareshian and Wachs
on the homology of Rees products of posets.

\end{abstract}

\section{Introduction and results}
\label{sec:intro}

Local $h$-polynomials were introduced by Stanley~\cite{Sta92} as a fundamental
tool in his theory of face enumeration for subdivisions of simplicial complexes.
Given a (finite, geometric) simplicial subdivision (triangulation) $\Gamma$ of
the abstract simplex $2^V$ on an $n$-element vertex set $V$, the local
$h$-polynomial $\ell_V (\Gamma, x)$ is defined by the formula

\begin{equation} \label{eq:deflocalh}
  \ell_V (\Gamma, x) \ = \sum_{F \subseteq V} \, (-1)^{n - |F|} \,
  h (\Gamma_F, x),
\end{equation}

\noindent
where $\Gamma_F$ is the restriction of $\Gamma$ to the face $F \in 2^V$ and
$h(\Delta, x)$ stands for the $h$-polynomial of the simplicial complex
$\Delta$; missing definitions can be found in Sections~\ref{sec:color} and
\ref{sec:sub} (the general notion of topological subdivision introduced in
\cite{Sta92} will not concern us here).

The importance of local $h$-polynomials stems from their appearance in the
locality formula \cite[Theorem~3.2]{Sta92}, which expresses the $h$-polynomial
of a simplicial subdivision of a pure simplicial complex $\Delta$ as a sum
of local contributions, one for each face of $\Delta$.

The polynomial $\ell_V (\Gamma, x)$ enjoys a number of pleasant properties,
interpretations and outstanding open problems. For instance, it was shown in
\cite[Sections~3--5]{Sta92} that $\ell_V (\Gamma, x)$ has nonnegative and
symmetric coefficients (with center of symmetry $n/2$) for every simplicial
subdivision $\Gamma$ of $2^V$ and unimodal coefficients for every regular
such subdivision. Moreover, it is conjectured \cite[Question~3.5]{Ath12}
\cite[Conjecture~5.4]{Sta92} that $\ell_V (\Gamma, x)$ has unimodal
coefficients for every simplicial subdivision $\Gamma$ of $2^V$ (although
this property fails for the more general class of quasi-geometric
subdivisions \cite[Example~3.4]{Ath12}) and that it is $\gamma$-nonnegative
for every flag such subdivision \cite[Conjecture~5.4]{Ath12}, meaning that
the integers $\xi_i (\Gamma)$ uniquely determined by

\begin{equation} \label{eq:deflocalg}
  \ell_V (\Gamma, x) \ = \sum_{i=0}^{\lfloor n/2 \rfloor} \, \xi_i (\Gamma)
  \, x^i (1+x)^{n-2i}
\end{equation}

\noindent
are nonnegative for all $0 \le i \le \lfloor n/2 \rfloor$. Combinatorial interpretations for the coefficients of $\ell_V (\Gamma, x)$ and the
integers $\xi_i (\Gamma)$ have been given for several interesting classes
of (flag) simplicial subdivisions of the simplex in \cite{AS12} \cite{Sav13}
\cite[Section~2]{Sta92}.

We now recall these interpretations for the (first, simplicial) barycentric
subdivision of the simplex, which will be generalized by the first main
result of this paper. We will denote by $\sd(2^V)$ the barycentric subdivision
of the simplex $2^V$, by $\mathfrak{S}_n$ the symmetric group of permutations
of $\{1, 2,\dots,n\}$ and by $\dD_n$ the set of derangements (permutations
without fixed points) in $\mathfrak{S}_n$. We recall that an excedance of $w \in \mathfrak{S}_n$ is an index $1 \le i \le n$ such that $w(i) > i$. A
descending run of $w \in \mathfrak{S}_n$ is a maximal string $\{a,
a+1,\dots,b\}$ of integers, such that $w(a) > w(a+1) > \cdots > w(b)$. The
following statement is a combination of results of \cite[Section~4]{AS12}
and \cite[Section~2]{Sta92}. The expression for the local $h$-polynomial
of $\sd(2^V)$ given in (\ref{eq:localgammabary}) may also be derived from an
unpublished result of Gessel \cite[Theorem~7.3]{SW10} (see
\cite[Remark~5.5]{SW10}), from Corollary 9 and Theorem 11 in \cite{SZ12} and, 
as will be explained in the sequel, from results of \cite{LSW12, SW09} on the 
homology of Rees products of posets.

\begin{theorem} \label{thm:bary} {\rm (\cite[Theorem~1.4]{AS12}
\cite[Proposition~2.4]{Sta92})}
Let $V$ be an $n$-element set. The local $h$-polynomial of the barycentric
subdivision of the $(n-1)$-dimensional simplex $2^V$ can be expressed as
  \begin{eqnarray} \label{eq:localbary}
    \ell_V (\sd(2^V), x) &=& \sum_{w \in \dD_n} x^{\ex(w)} \\
  & & \nonumber \\
  &=& \label{eq:localgammabary}
      \sum_{i=1}^{\lfloor n/2 \rfloor} \xi_{n, i} \, x^i (1+x)^{n-2i},
  \end{eqnarray}
where $\dD_n$ is the set of derangements in $\mathfrak{S}_n$, $\ex(w)$ is
the number of excedances of $w \in \mathfrak{S}_n$ and $\xi_{n,i}$ stands
for the number of permutations $w \in \mathfrak{S}_n$ with $i$ descending
runs and no descending run of size one.
\end{theorem}

  \begin{figure}
  \epsfysize = 1.6 in \centerline{\epsffile{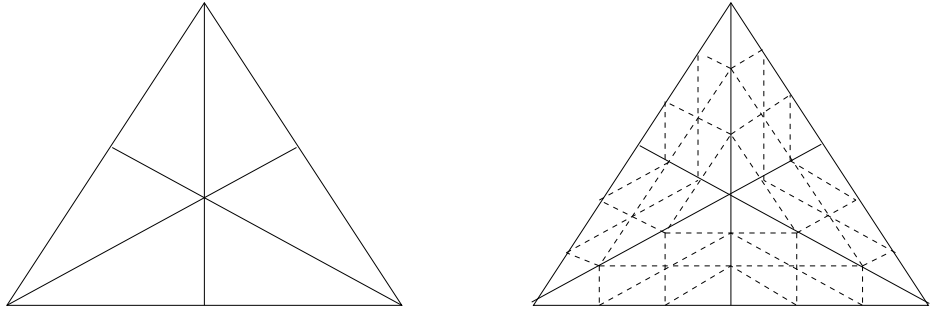}}
  \caption{The barycentric subdivision of the 2-simplex and its third
  edgewise subdivision}
  \label{fig:KK33}
  \end{figure}

Our generalization of Theorem~\ref{thm:bary} concerns the $r$th edgewise
subdivision $\sd(2^V)^{\langle r \rangle}$ of the barycentric subdivision
$\sd(2^V)$ and the flag excedance statistic on the wreath product $\ZZ_r
\wr \mathfrak{S}_n$, where $\ZZ_r$ is the cyclic group of order $r$. The
$r$th edgewise subdivision $\Delta^{\langle r \rangle}$ of a simplicial
complex $\Delta$ is a standard way to subdivide a simplicial complex
$\Delta$ so that each face $F \in \Delta$ is subdivided into $r^{\dim(F)}$
faces of the same dimension. This construction has appeared is several
mathematical contexts; see, for instance, \cite{BrW09, BR05, EG00, Gr89}.
Figure~\ref{fig:KK33} shows the barycentric subdivision of the 2-dimensional
simplex and its third edgewise subdivision. The elements of $\ZZ_r \wr
\mathfrak{S}_n$ will be thought of as $r$-colored permutations, meaning
permutations in $\mathfrak{S}_n$ with each coordinate in their one--line
notation colored with an element of the set $\{0, 1,\dots,r-1\}$. The
flag excedance of $w \in \ZZ_r \wr \mathfrak{S}_n$ is defined as

\begin{equation} \label{eq:def-fexc}
  \fex(w) \ = \ r \cdot \ex_A (w) \, + \, \csum(w),
\end{equation}

\noindent
where $\ex_A (w)$ is the number of excedances of $w$ at coordinates of zero
color and $\csum(w)$ is the sum of the colors of all coordinates of $w$; more
explanation appears in Section~\ref{sec:color}. The flag excedance statistic
was introduced (although not with this name) by Bagno and Garber \cite{BG06}
and was further studied in \cite{FH11} and, in the special case $r=2$, in
\cite{FH09, Mo13}. We will denote by $\dD^r_n$ the set of derangements
(elements without fixed points of zero color) in $\ZZ_r \wr \mathfrak{S}_n$.
We will call an element $w \in \ZZ_r \wr \mathfrak{S}_n$ balanced if $\csum(w)$
(equivalently, $\fex(w)$) is divisible by $r$. Decreasing runs for elements
of $\ZZ_r \wr \mathfrak{S}_n$ are defined as for those of $\mathfrak{S}_n$,
after the $r$-colored integers have been totally ordered in a standard
lexicographic way; see Section~\ref{sec:color}. The following statement
provides a new context in which the flag excedance statistic naturally
arises; it reduces to Theorem~\ref{thm:bary} for $r=1$.

\begin{theorem} \label{thm:edgebary}
Let $V$ be an $n$-element set and let $r$ be a positive integer. The local
$h$-polynomial of the $r$th edgewise subdivision of the barycentric
subdivision of the $(n-1)$-dimensional simplex $2^V$ can be expressed as
  \begin{eqnarray} \label{eq:localedgebary}
    \ell_V (\sd(2^V)^{\langle r \rangle}, x) &=& \sum_{w \in (\dD^r_n)^b}
    x^{\fex (w) / r} \\
  & & \nonumber \\
  &=& \label{eq:localgammaedgebary}
      \sum_{i=1}^{\lfloor n/2 \rfloor} \xi^+_{n, r, i} \, x^i (1+x)^{n-2i},
  \end{eqnarray}
where $(\dD^r_n)^b$ is the set of balanced derangements in $\ZZ_r \wr
\mathfrak{S}_n$, $\fex (w)$ is the flag excedance of $w \in \ZZ_r \wr
\mathfrak{S}_n$ and $\xi^+_{n, r, i}$ stands for the number of elements of
$\ZZ_r \wr \mathfrak{S}_n$ with $i$ descending runs, no descending run of
size one and last coordinate of zero color.
\end{theorem}

Theorem~\ref{thm:edgebary} will be derived from a related result on the
derangement polynomial for $\ZZ_r \wr \mathfrak{S}_n$, studied by Chow
and Mansour \cite{CM10}. This polynomial, denoted here by $d^r_n (x)$, is
defined as the generating polynomial for the excedance statistic, in the
sense of Steingr\'imsson~\cite{Stei92, Stei94}, on the set $\dD^r_n$ of
derangements in $\ZZ_r \wr \mathfrak{S}_n$. It specializes to the
right-hand side of (\ref{eq:localbary}), the derangement polynomial for
$\mathfrak{S}_n$, introduced by Brenti~\cite{Bre90}, for $r=1$ and to its
type $B$ analogue introduced by Chen, Tang and Zhao~\cite{CTZ09} and,
independently (in a variant form), by Chow~\cite{Ch09}, for $r=2$. The
following statement is the second main result of this paper.

\begin{theorem} \label{thm:dnr}
For all positive integers $n, r$
  \begin{equation} \label{eq:dnrdecompose}
    d^r_n (x) \ = \ \sum_{i=0}^{\lfloor n/2 \rfloor} \xi^+_{n, r, i} \, x^i
     (1+x)^{n-2i} \ \, + \sum_{i=0}^{\lfloor (n+1)/2 \rfloor} \xi^-_{n, r, i}
       \, x^i (1+x)^{n+1-2i},
  \end{equation}
where $\xi^+_{n, r, i}$ is as in Theorem~\ref{thm:edgebary} and $\xi^-_{n,
r, i}$ stands for the number of elements of $\ZZ_r \wr \mathfrak{S}_n$ with
$i$ descending runs, no descending run of size one other than $\{n\}$
and last coordinate of nonzero color.
\end{theorem}

The previous theorem implies that $d^r_n (x)$ can be written as a sum of
two symmetric and unimodal (in fact, $\gamma$-nonnegative) polynomials whose
centers of symmetry differ by a half. Moreover, one of these polynomials is
equal to the local $h$-polynomial of a suitable simplicial subdivision of
the $(n-1)$-dimensional simplex. These facts were first observed in the
special case $r=2$ in joint work of the author with Savvidou
\cite[Chapter~3]{Sav13} using methods different from those in this paper
(the two main results of this paper, however, are new even for $r=2$).
Although a direct combinatorial proof should be possible,
Theorem~\ref{thm:dnr} will be derived from \cite[Corollary~3.8]{LSW12} by
exploiting a connection of $d^r_n (x)$ to the homology of Rees products of
posets, previously noticed only in the special case $r=1$.

\medskip
This paper is structured as follows. Section~\ref{sec:color}
includes background material on the combinatorics of colored
permutations. Some preliminary results on the various excedance
statistics on derangements in $\ZZ_r \wr \mathfrak{S}_n$ are also
included. Section~\ref{sec:rees} reviews the necessary background
on Rees product homology and proves Theorem~\ref{thm:dnr}.
Theorem~\ref{thm:edgebary} is proven in Section~\ref{sec:proof}
after the relevant definitions on simplicial subdivisions and
local $h$-polynomials have been explained in
Section~\ref{sec:sub}. Section~\ref{sec:rem} concludes with some
corollaries, remarks and open problems.

\section{Colored permutation statistics}
\label{sec:color}

This section recalls basic definitions and useful facts about colored
permutations and the various notions of descent and excedance for them. It
also includes some new results and formulas, in preparation for the proofs
of Theorems~\ref{thm:edgebary} and \ref{thm:dnr} in the following sections.

\subsection{Colored permutations}
\label{subsec:cperm}

The symmetric group of permutations of $\{1, 2,\dots,n\}$ will be denoted
by $\mathfrak{S}_n$. The elements of the cyclic group $\ZZ_r$ of order $r$
will be represented by those of $\{0, 1,\dots,r-1\}$ and will be thought
of as colors.

The wreath product $\ZZ_r \wr \mathfrak{S}_n$ consists of all $r$-colored
permutations of the form $\sigma \times \rzz$, where $\sigma = (\sigma(1),
\sigma(2),\dots,\sigma(n)) \in \mathfrak{S}_n$ and $\rzz = (z_1,
z_2,\dots,z_n) \in \{0, 1,\dots,r-1\}^n$. The number $z_i$ will be thought
of as the color assigned to $\sigma(i)$. The product in the group $\ZZ_r
\wr \mathfrak{S}_n$ is given by the rule $(\sigma \times \rzz) (\tau \times
\rww) = \sigma \tau \times (\rww + \tau(\rzz))$, where the composition
$\sigma \tau = \sigma \circ \tau$ is evaluated from right to left, $\tau
(\rzz) = (z_{\tau(1)}, z_{\tau(2)},\dots,z_{\tau(n)})$ and the addition is
coordinatewise modulo $r$. The inverse of $\sigma \times \rzz \in \ZZ_r \wr
\mathfrak{S}_n$ is the element $\sigma^{-1} \times \rww$, where
$\sigma^{-1}$ is the inverse of $\sigma$ in $\mathfrak{S}_n$ and $\rww =
(w_1, w_2,\dots,w_n)$ is such that

  $$ w_i \ = \ \cases{ 0, & if \ $z_{\sigma^{-1}(i)} = 0$, \cr
      r - z_{\sigma^{-1}(i)}, & otherwise. } $$

\noindent
A \emph{fixed point} of $\sigma \times \rzz \in \ZZ_r \wr \mathfrak{S}_n$
is an index $i \in \{1, 2,\dots,n\}$ such that $\sigma(i) = i$ and $z_i =
0$. Elements of $\ZZ_r \wr \mathfrak{S}_n$ without fixed points are called
\emph{derangements}. The set of derangements in $\ZZ_r \wr \mathfrak{S}_n$
will be denoted by $\dD^r_n$.

\subsection{Statistics}
\label{subsec:stat}

Throughout this section $w = \sigma \times \rzz \in \ZZ_r \wr \mathfrak{S}_n$
will be a colored permutation as in Section~\ref{subsec:cperm}, where $\sigma
= (\sigma(1), \sigma(2),\dots,\sigma(n)) \in \mathfrak{S}_n$ and $\rzz = (z_1,
z_2,\dots,z_n) \in \{0, 1,\dots,r-1\}^n$.

We will first recall the notions of descent and excedance for elements of
$\ZZ_r \wr \mathfrak{S}_n$, introduced in \cite[Section~2]{Stei92}
\cite[Section~3.1]{Stei94}. A \emph{descent} of $w$ is an index $i \in \{1,
2,\dots,n\}$ such that either $z_i > z_{i+1}$, or $z_i = z_{i+1}$ and $\sigma(i)
> \sigma(i+1)$, where $\sigma(n+1) = n+1$ and $z_{n+1} = 0$. Thus, in
particular, $n$ is a descent of $w$ if and only if $\sigma(n)$ has nonzero
color. A \emph{descending run} of $w$ is a maximal string $\{a, a+1,\dots,b\}$
of integers such that $i$ is a descent of $w$ for every $a \le i \le b-1$.
An \emph{excedance} of $w$ is an index $i \in \{1, 2,\dots,n\}$ such that
either $\sigma(i) > i$, or $\sigma(i) = i$ and $z_i > 0$. The number of
descents and excedances, respectively, of $w$ will be denoted by $\des (w)$
and $\ex(w)$. The equidistribution of the descent and excedance statistics
\cite[Theorem~3.15]{Stei92} \cite[Theorem~15]{Stei94} on $\ZZ_r \wr
\mathfrak{S}_n$ allows one to define the Eulerian polynomial for $\ZZ_r \wr
\mathfrak{S}_n$ as

  \begin{equation} \label{eq:stei}
    A^r_n (x) \ := \sum_{w \in \ZZ_r \wr \mathfrak{S}_n} x^{\des (w)} \ =
      \sum_{w \in \ZZ_r \wr \mathfrak{S}_n} x^{\ex (w)}.
  \end{equation}

\noindent
For $r=1$, this polynomial specializes to the classicial Eulerian polynomial

  $$ A_n (x) \ := \sum_{w \in \mathfrak{S}_n} x^{\des (w)} \ = \sum_{w \in
     \mathfrak{S}_n} x^{\ex (w)}, $$

\noindent
although this definition differs from the standard one \cite[p.~32]{StaEC1}
by a power of $x$.

The \emph{flag descent} statistic on $\ZZ_r \wr \mathfrak{S}_n$ was defined
first in the special case $r=2$ by Adin, Brenti and
Roichman~\cite[Section~4]{ABR01} and later by Bagno and
Biagoli~\cite[Section~3]{BB07} for general $r \ge 1$. We find it convenient
to define it here as

  \begin{equation} \label{eq:def-fdes}
    \fdes^*(w) \ = \ \cases{ r \cdot \des (w), & if \ $z_n = 0$, \cr
                           r \cdot \des (w) - r + z_n, & otherwise }
  \end{equation}

\noindent
for $w \in \ZZ_r \wr \mathfrak{S}_n$ where, as explained earlier, $z_n \in \{0, 1,\dots,r-1\}$ is the color of the last coordinate 
$\sigma(n)$ of $w$ (this definition differs by a minor twist from the one given in \cite{ABR01, BB07}; see the proof of 
Proposition~\ref{prop:fdesfexc}). The \emph{flag excedance} number $\fex(w)$ of $w \in \ZZ_r \wr \mathfrak{S}_n$ is defined by
(\ref{eq:def-fexc}), where $\ex_A (w)$ is the number of indices $i \in \{1, 2,\dots,n\}$ such that $\sigma(i) > i$ and $z_i = 0$ and $\csum(w) = z_1 + z_2 + \cdots + z_n$, the sum being computed in $\ZZ$. We will call $w$ \emph{balanced} if $\csum(w)$ is an integer multiple of $r$.

\begin{example} \rm
Suppose that $w = \sigma \times \rzz$ is the colored permutation in $\ZZ_3
\wr \mathfrak{S}_6$ with $\sigma = (2, 5, 1, 4, 6, 3)$ and $\rzz = (1, 0,
0, 2, 0, 1)$. Then $w$ is a derangement with descents 1, 2, 4 and 6,
excedances 1, 2, 4 and 5 and type $A$ excedances 2 and 5, so that $\des(w)
= \ex(w) = 4$, $\ex_A(w) = 2$ and $\fdes^*(w) = 10$. Moreover, we have
$\csum(w) = 4$ and hence $\fex(w) = 10$.
\end{example}

The following statement combines results from \cite{ABR01, BB07, FH11}
(where \cite{ABR01} contributes to the special case $r=2$).

\begin{proposition} \label{prop:fdesfexc} {\rm (\cite{ABR01, BB07, FH11})}
For all positive integers $n, r$

  \begin{equation} \label{eq:fdesfexc}
    \sum_{w \in \ZZ_r \wr \mathfrak{S}_n} x^{\fdes^* (w)} \ = \sum_{w \in
      \ZZ_r \wr \mathfrak{S}_n} x^{\fex (w)} \ = \
        (1 + x + x^2 + \cdots + x^{r-1})^n \, A_n (x).
  \end{equation}
\end{proposition}

\noindent
\emph{Proof.} Given $w = \sigma \times \rzz \in \ZZ_r \wr \mathfrak{S}_n$,
as in the beginning of this section, we set $\tilde{w} = \tilde{\sigma}
\times \tilde{\rzz} \in \ZZ_r \wr \mathfrak{S}_n$ where $\tilde{\sigma}
= (n - \sigma(n) + 1, n - \sigma(n-1) + 1,\dots,n - \sigma(1) + 1) \in
\mathfrak{S}_n$ and $\tilde{\rzz} = (z_n, z_{n-1},\dots,z_1)$.
The flag descent of $w$ is defined in \cite[Section~3]{BB07} as $\fdes(w)
= r \cdot \des^*(w) + z_1$, where $\des^*(w)$ is the number of indices
$i \in \{1, 2,\dots,n-1\}$ such that either $z_i < z_{i+1}$, or $z_i =
z_{i+1}$ and $\sigma(i) > \sigma(i+1)$. We note that this happens if and
only if $n-i$ is a descent of $\tilde{w}$ and conclude that

 $$ \fdes(w) \ = \ r \cdot \des_A (\tilde{w}) \, + \, z_n(\tilde{w}) \ =
    \ \fdes^*(\tilde{w}), $$

\noindent
where $\des_A (\tilde{w})$ is the number of descents of $\tilde{w}$ other
than $n$. Since the map which sends $w$ to $\tilde{w}$ is a bijection from
$\ZZ_r \wr \mathfrak{S}_n$ to itself, this implies that

  $$ \sum_{w \in \ZZ_r \wr \mathfrak{S}_n} x^{\fdes^* (w)} \ = \sum_{w \in
      \ZZ_r \wr \mathfrak{S}_n} x^{\fdes (w)}. $$

\noindent
Given that, the equality between the leftmost and rightmost expressions
in (\ref{eq:fdesfexc}) follows by setting $q=1$ in \cite[Theorem A.1]{BB07}.
The first equality in (\ref{eq:fdesfexc}) follows from
\cite[Theorem~1.4]{FH11}.
\qed

\bigskip
We now consider the generating polynomials for the excedance and flag
excedance statistics on the set $\dD^r_n$ of derangements in $\ZZ_r \wr
\mathfrak{S}_n$. The generating polynomial

  \begin{equation} \label{eq:defdnr}
    d_n^r (x) \ := \, \sum_{w \in \dD^r_n} x^{\ex (w)}
  \end{equation}

\noindent for the excedance statistic was studied by Chow and
Mansour \cite[Section~3]{CM10}. As already mentioned in the
introduction, $d_n^r (x)$ generalizes both the right-hand side of
(\ref{eq:localbary}) and (as can be inferred, for instance, from
the discussion in \cite[p.~428]{BB07}) the type $B$ derangement
polynomial, introduced and studied in \cite{CTZ09, Ch09}. For $x
\in \RR$, we will denote by $\lceil x \rceil$ the smallest integer
which is not strictly less than $x$.

\begin{proposition} \label{prop:dnr1}
For all positive integers $n, r$
  \begin{equation} \label{eq:dnr1}
     d_n^r (x) \ = \, \sum_{w \in \dD^r_n} x^{\lceil \frac{\fex(w)}{r} \rceil}.
  \end{equation}
\end{proposition}

\noindent
\emph{Proof.} We define the linear operator $\widetilde{{\rm E}}_r: \RR[x]
\to \RR[x]$ by setting $\widetilde{{\rm E}}_r (x^m) = x^{\lceil m/r \rceil}$
for $m \in \NN$. The defining equation (\ref{eq:def-fdes}) shows that
$\widetilde{{\rm E}}_r (x^{\fdes^*(w)}) = x^{\des(w)}$ for every $w \in \ZZ_r
\wr \mathfrak{S}_n$. Applying operator $\widetilde{{\rm E}}_r$ to the first
equality in (\ref{eq:fdesfexc}) we get

  $$ A^r_n (x) \ = \sum_{w \in \ZZ_r \wr \mathfrak{S}_n} x^{\des (w)} \ =
     \sum_{w \in \ZZ_r \wr \mathfrak{S}_n} x^{\lceil \frac{\fex(w)}{r} \rceil}.
  $$

\noindent
Because of (\ref{eq:stei}), the previous equalities can be rewritten as

  $$ A^r_n (x) \ = \sum_{w \in \ZZ_r \wr \mathfrak{S}_n} x^{\ex (w)} \ =
     \sum_{w \in \ZZ_r \wr \mathfrak{S}_n} x^{\lceil \frac{\fex(w)}{r} \rceil}.
  $$

\noindent
Since adding or removing fixed points from $w \in \ZZ_r \wr \mathfrak{S}_n$
does not affect $\ex(w)$ or $\fex(w)$, the previous equalities, the defining
equation (\ref{eq:defdnr}) and two easy applications of the principle of
Inclusion-Exclusion show that
  \begin{eqnarray*}
    d_n^r (x) &=& \sum_{k=0}^n \, (-1)^{n-k} {n \choose k}
                  \sum_{w \in \ZZ_r \wr \mathfrak{S}_k} x^{\ex (w)} \\
  & & \\
  &=& \sum_{k=0}^n \, (-1)^{n-k} {n \choose k} \sum_{w \in \ZZ_r \wr
      \mathfrak{S}_k} x^{\lceil \frac{\fex(w)}{r} \rceil} \\
  & & \\
  &=& \sum_{w \in \dD^r_n} x^{\lceil \frac{\fex(w)}{r} \rceil}
  \end{eqnarray*}

\noindent
and the proof follows.
\qed

\begin{proposition} \label{prop:dnr2}
For all positive integers $n, r$

  \begin{equation} \label{eq:dnr2}
    x^n d_n^r (1/x) \ = \ \sum_{k=0}^n \, (-1)^{n-k} {n \choose k} r^k A_k(x),
  \end{equation}
where $A_0 (x) = 1$.
\end{proposition}

\noindent
\emph{Proof.} It was shown in \cite[Theorem~5~(iv)]{CM10} that

  $$  \sum_{n \ge 0} \, d^r_n (x) \, \frac{t^n}{n!} \ = \ \frac{(1-x) e^{(r-1)xt}}
    {e^{rxt} - xe^{rt}}. $$

\noindent
Replacing $x$ by $1/x$ and $t$ by $xt$ we get

  \begin{equation} \label{eq:dnrexp}
    \sum_{n \ge 0} \, x^n d^r_n (1/x) \, \frac{t^n}{n!} \ = \ \frac{(1-x) e^{(r-1)t}}
    {e^{rxt} - xe^{rt}}.
  \end{equation}

Similarly, denoting by $S_n (x)$ the right-hand side of (\ref{eq:dnr2}) and using
the well-known expression for the exponential generating function for the Eulerian
polynomials $A_n (x)$ \cite[Proposition~1.4.5]{StaEC1} we get

  $$  \sum_{n \ge 0} \, S_n (x) \, \frac{t^n}{n!} \ = \ e^{-t} \, \cdot \,
      \sum_{n \ge 0} \, A_n (x) \, \frac{(rt)^n}{n!}  \ = \ \frac{(1-x) e^{(r-1)t}}
    {e^{rxt} - xe^{rt}}. $$

\noindent
Comparing with (\ref{eq:dnrexp}) shows that $x^n d_n^r (1/x) = S_n (x)$ for every
$n \in \NN$ and the proof follows.
\qed

\bigskip
The following observation on the generating polynomial for the flag excedance
statistic on $\dD^r_n$ specializes to \cite[Proposition~3.5]{Mo13} for $r=2$.

\begin{proposition} \label{prop:dnrflag}
For all positive integers $n, r$, the polynomial
  \begin{equation} \label{eq:defdnrflag}
    f^r_n (x) \ := \, \sum_{w \in \dD^r_n} x^{\fex (w)}
  \end{equation}
is symmetric with center of symmetry $rn/2$, i.e., we have $x^{rn} f^r_n
(1/x) = f^r_n (x)$. Moreover, $f^r_n (x)$ is monic of degree $rn-1$.
\end{proposition}

\noindent \emph{Proof.} Suppose that $w \in \dD^r_n$ is a colored
permutation with $k$ coordinates of zero color. From the
description of $w^{-1}$ given in Section~\ref{subsec:cperm} we get
$\ex_A (w^{-1}) = k - \ex_A(w)$ and $\csum(w^{-1}) = r(n-k) -
\csum(w)$. These equalities imply that $\fex (w^{-1}) = rn -
\fex(w)$. Since the map $\varphi: \dD^r_n \to \dD^r_n$ defined by
$\varphi(w) = w^{-1}$ for $w \in \dD^r_n$ is a well-defined
bijection, the first claim in the statement of the proposition
follows. The second can be left to the reader. \qed

\section{Rees product homology and proof of Theorem~\ref{thm:dnr}}
\label{sec:rees}

This section derives Theorem~\ref{thm:dnr} from results of Shareshian and
Wachs~\cite{SW09} and Linusson, Shareshian and Wachs~\cite{LSW12} on the
homology of Rees products of posets. Some of the background needed to
understand and apply these results will first be explained. Throughout
this section we will assume familiarity with the basics of finite partially
ordered sets \cite[Chapter~3]{StaEC1}. We will also write $[a, b] = \{a,
a+1,\dots,b\}$ for integers $a, b$ with $a \le b$ and set $[n] := [1, n]$.

\subsection{Shellability and Rees products of posets}
\label{subsec:shellrees}

Given a finite poset $P$ with a minimum element $\hat{0}$ and maximum
element $\hat{1}$, we will denote by $\bar{P}$ the poset which is obtained
from $P$ by removing $\hat{0}$ and $\hat{1}$. Similarly, given any finite
poset $Q$, we will denote by $\hat{Q}$ the poset which is obtained from $Q$
by adding elements $\hat{0}$ and $\hat{1}$ so that $\hat{0} < x < \hat{1}$
holds for every $x \in Q$. We will then write $\mu (Q)$ for the value
$\mu_{\hat{Q}} (\hat{0}, \hat{1})$ of the M\"obius function
\cite[Section~3.7]{StaEC1} of $\hat{Q}$ between $\hat{0}$ and $\hat{1}$.
We will denote by $Q_{\le x}$ the principal order ideal $\{y \in Q: y \le
x\}$ of $Q$ generated by $x \in Q$.

We now briefly recall the definition of EL-shellability and refer the reader
to \cite{Bj95} \cite[Lecture~3]{Wa07} for more information. Suppose $P$ is
a finite graded poset of rank $n+1$, with minimum element $\hat{0}$ and
maximum element $\hat{1}$. Consider a map $\lambda: \eE(P) \to \Lambda$,
where $\eE(P)$ is the set of cover relations of $P$ and $\Lambda$ is a set
equipped with a total order $\preceq$. To any unrefinable chain $c: x_0 <
x_1 < \cdots < x_r$ of elements of $P$ one can associate the sequence
$\lambda(c) = (\lambda(x_0, x_1), \lambda(x_1, x_2),\dots,\lambda(x_{r-1},
x_r))$, called here the \emph{label} of $c$. The chain $c$ is said to be
\emph{rising} or
\emph{falling} if the entries of $\lambda(c)$ weakly increase or strictly
decrease, respectively, in the total order $\preceq$. The map $\lambda$ is
called an \emph{EL-labeling} if the following hold for every non-singleton
interval $[u, v]$ in $P$: (a) there is a unique rising maximal chain in
$[u, v]$; and (b) the label of this chain is lexicographically smallest
among all labels of maximal chains in $[u, v]$. The poset $P$ is called
\emph{EL-shellable} if it has an EL-labeling for some totally ordered set
$\Lambda$. Assuming $P$ is EL-shellable and given $S \subseteq [n]$, we
will denote by $\beta_P (S)$ the number of maximal chains $c: \hat{0} =
x_0 < x_1 < \cdots < x_{n+1} = \hat{1}$ in $P$ such that for $1 \le i \le
n$ we have $\lambda(x_{i-1}, x_i) \succ \lambda(x_i, x_{i+1})
\Leftrightarrow i \in S$. The numbers $\beta_P (S)$ are independent of
the EL-labeling $\lambda$.

Given two finite graded posets $P$ and $Q$ with rank functions $\rho_P$ and
$\rho_Q$, respectively, the \emph{Rees product} of $P$ and $Q$ is defined
in \cite{BjW05} as $P \ast Q = \{ (p, q) \in P \times Q: \rho_P (p) \ge
\rho_Q (q) \}$, with partial order defined by setting $(p_1, q_1) \le (p_2,
q_2)$ if all of the following conditions are satisfied:

\begin{itemize}
\itemsep=0pt
\item[$\bullet$]
$p_1 \le p_2$ holds in $P$,

\item[$\bullet$]
$q_1 \le q_2$ holds in $Q$ and

\item[$\bullet$]
$\rho_P (p_2) - \rho_P (p_1) \ge \rho_Q (q_2) - \rho_Q (q_1)$.
\end{itemize}

\noindent
Equivalently, $(p_1, q_1)$ is covered by $(p_2, q_2)$ in $P \ast Q$ if and
only if (a) $p_1$ is covered by $p_2$ in $P$; and (b) either $q_1 = q_2$, or
$q_1$ is covered by $q_2$ in $Q$. As a consequence of the definition, the
principal order ideal of $P \ast Q$ generated by $x = (p, q) \in P \ast Q$
satisfies
  \begin{equation} \label{eq:reesideal}
    (P \ast Q)_{\le x} \ = \ (P_{\le p} \ast Q_{\le q})_{\le x}.
  \end{equation}

For positive integers $x, n$ we will denote by $T_{x, n}$ the poset whose
Hasse diagram is a complete $x$-ary tree of height $n-1$, with root at the
bottom (so that $T_{x, n}$ is an $n$-element chain for $x=1$). The following
theorem is a restatement of the first part of \cite[Corollary~3.8]{LSW12}.

\begin{theorem} \label{thm:LSW} {\rm (\cite{LSW12})}
For every EL-shellable poset $P$ of rank $n+1$ and every positive integer $x$
we have
  \begin{eqnarray} \label{eq:LSWsum}
    |\mu (\bar{P} \ast T_{x, n})| &=& \sum_{S \in \pP_\stab ([2, n-1])}
    \beta_P ([n] \sm S) \, x^{|S|} \, (1+x)^{n-1-2|S|} \ \ + \\
  & & \nonumber \\
  & & \sum_{S \in \pP_\stab ([2, n-2])} \beta_P ([n-1] \sm S) \, x^{|S|+1} \,
      (1+x)^{n-2-2|S|}, \nonumber
  \end{eqnarray}

\noindent
where $\pP_\stab (\Theta)$ denotes the set of all subsets of $\Theta \subseteq
\ZZ$ which do not contain two consecutive integers.
\end{theorem}

\subsection{Proof of Theorem~\ref{thm:dnr}}
\label{subsec:rees}

We consider the subsets $\Omega$ of $[n] \times [0, r-1]$ for which for every
$i \in [n]$ there is at most one $j \in [0, r-1]$ such that $(i, j) \in \Omega$.
We will denote by $B^r_n$ the set of all nonempty such subsets, partially
ordered by inclusion. Thus, for $r=1$ the poset $B^r_n$ is isomorphic to the
poset obtained from the Boolean lattice of subsets of $[n]$ by removing its
minimum element (empty set). For $r=2$ this poset was considered in
\cite[Section~6]{SW09}.

The following statement provides the key connection of the
derangement polynomial $d_n^r (x)$ to the Rees product
construction. The special case $r=1$ is equivalent to the special
case $q=1$ of \cite[Equation~(1.3)]{LSW12}.

\begin{proposition} \label{prop:reesdnr}
For all positive integers $n, r, x$ we have
  \begin{equation} \label{eq:reesdnr}
     |\mu (B^r_n \ast T_{x, n})| \ = \ x^n d_n^r (1/x).
  \end{equation}
\end{proposition}

\noindent
\emph{Proof.} Setting $Q = B^r_n \ast T_{x, n}$ and using the definition
\cite[Section~3.7]{StaEC1} of the M\"obius function we get

\begin{equation} \label{eq:mudef}
  \mu (B^r_n \ast T_{x, n}) \ = \ \mu(Q) \ = \ \mu_{\hat{Q}} (\hat{0}, \hat{1})
     \ = \ - \sum_{y \in Q \cup
\{ \hat{0} \}} \mu_{\hat{Q}} (\hat{0}, y).
\end{equation}

\noindent
Suppose that $y = (a, b) \in Q$, where $a \in B^r_n$ has rank $k-1$, $b \in T_{x,
n}$ has rank $i$ and $k-1, i \in [0, n-1]$ with $i \le k-1$. Since the
principal order ideal of $B^r_n$ generated by $a$ is isomorphic to $B^1_k$ and
that of $T_{x, n}$ generated by $b$ is an $(i+1)$-element chain,
Equation~(\ref{eq:reesideal}) and \cite[Theorem~1.2]{SW09} imply that $(-1)^k
\mu_{\hat{Q}} (\hat{0}, y)$ is equal to the (Eulerian) number of elements of
$\mathfrak{S}_k$ with $i$ descents. Denoting this number by $a_{ki}$ and noting
that there are exactly $r^k {n \choose k}$ elements of $B^r_n$ of rank $k-1$
and exactly $x^i$ elements of $T_{x, n}$ of rank $i$, we conclude from
(\ref{eq:mudef}) that
  \begin{eqnarray*}
    \mu (Q) &=& -1 \, - \, \sum_{i=0}^{n-1} \, \sum_{k=i+1}^n r^k {n \choose
                  k} x^i \, (-1)^k a_{ki} \\
    & & \\
    &=& -1 \, - \, \sum_{k=1}^n \, (-1)^k {n \choose k} r^k \, \sum_{i=0}^{k-1}
        a_{ki} x^i \\
    & & \\
    &=& \sum_{k=0}^n \, (-1)^{k+1} {n \choose k} r^k A_k (x)
  \end{eqnarray*}

\noindent
and the result follows from Proposition~\ref{prop:dnr2}.
\qed

\bigskip
\noindent
\emph{Proof of Theorem \ref{thm:dnr}.} Since both sides of
(\ref{eq:dnrdecompose}) are polynomials in $x$, we may assume that 
$x$ is a positive integer. We wish to apply Theorem~\ref{thm:LSW} to the poset $P := \hat{B^r_n}$.
Clearly, $P$ is graded of rank $n+1$. We set $\Lambda = ([n] \times [0, r-1]) \cup \{\delta\}$, where $\delta$ is a symbol not in $[n] \times [0, r-1]$, and totally order $\Lambda$ as
\begin{center}
\begin{tabular}{c}
$(1, 0) \prec (2, 0) \prec \cdots \prec (n, 0) \prec \delta \prec
     (1, 1) \prec (2, 1) \prec \cdots \prec (n, 1) \prec \cdots \prec
     (1, r-1) \prec$ \\ 
$(2, r-1) \prec \cdots \prec (n, r-1). $
\end{tabular} 
\end{center}
For a cover relation $(x, y) \in \eE(P)$ we define
$\lambda(x, y)$ as the unique element of the set $y \sm x$, if $y
\in P \sm \{ \hat{1} \}$, and set $\lambda(x, y) = \delta$ if $y =
\hat{1}$. We leave to the reader to verify (by standard arguments; see, for instance, \cite[Section~3.2.1]{Wa07}) that $\lambda: \eE(P)
\to \Lambda$ is an EL-labeling and that the map 
which assigns to a maximal chain $c$ of $P$ the label $\lambda(c)$, with its last entry $\delta$ removed, induces a bijection from the set of maximal chains of $P$ to $\ZZ_r \wr \mathfrak{S}_n$. Moreover,
given such a chain $c: \hat{0} = x_0 < x_1 < \cdots < x_{n+1} =
\hat{1}$ with corresponding colored permutation $w$, an index $i
\in [n]$ satisfies $\lambda(x_{i-1}, x_i) \succ \lambda(x_i,
x_{i+1})$ if and only if $i$ is a descent of $w$, in the sense of
Section~\ref{subsec:stat}. As a result, for $S \subseteq [n]$ the
number $\beta_P(S)$ is equal to the number of elements of $\ZZ_r
\wr \mathfrak{S}_n$ whose set of descents is equal to $S$.

Combining Proposition~\ref{prop:reesdnr} with Theorem~\ref{thm:LSW} and the
previous discussion (and remembering to replace $x$ by $1/x$ in the right-hand
side of (\ref{eq:LSWsum})), we conclude that (\ref{eq:dnrdecompose}) holds for
some nonnegative integers $\xi^+_{n, r, i}$ and $\xi^-_{n, r, i}$. Moreover,
denoting by $\Asc(w)$ the set of elements of $[n]$ which are not descents of $w \in \ZZ_r \wr \mathfrak{S}_n$, we have

  \begin{eqnarray*}
    \xi^+_{n, r, i} &=& \sum \, \beta_P ([n-1] \sm S) \\
    &=& \sum \ \# \ \{ w \in \ZZ_r \wr \mathfrak{S}_n: \, \Asc(w) = S \cup
        \{n\} \},
  \end{eqnarray*}

\medskip
\noindent
where $S$ runs through all $(i-1)$-element sets in $\pP_\stab ([2, n-2])$ and

  \begin{eqnarray*}
    \xi^-_{n, r, i} &=& \sum \, \beta_P ([n] \sm S) \\
    &=& \sum \ \# \ \{ w \in \ZZ_r \wr \mathfrak{S}_n: \, \Asc(w) = S \},
  \end{eqnarray*}

\medskip
\noindent
where $S$ runs through all $(i-1)$-element sets in $\pP_\stab ([2, n-1])$. These
are exactly the combinatorial interpretations claimed for $\xi^+_{n, r, i}$ and $\xi^-_{n, r, i}$ in Theorem \ref{thm:dnr}.
\qed

\section{Subdivisions}
\label{sec:sub}

This section briefly reviews the background on simplicial subdivisions needed
to understand Theorem~\ref{thm:edgebary} and its proof. Familiarity with basic
notions on simplicial complexes, such as the correspondence between abstract and
geometric simplicial complexes, will be assumed. For more information on this
topic and any undefined terminology, the reader is referred to \cite{Bj95,
Mat03, StaCCA}. All simplicial complexes considered here will be finite. We
will denote by $|S|$ the cardinality, and by $2^S$ the set of all subsets, of
a finite set $S$.

Consider two geometric simplicial complexes $\Sigma'$ and $\Sigma$ in some
Euclidean space $\RR^m$ (so the elements of $\Sigma'$ and $\Sigma$ are
geometric simplices in $\RR^m$), with corresponding abstract simplicial
complexes $\Delta'$ and $\Delta$. We will say that $\Sigma'$ is a
\emph{simplicial subdivision} of $\Sigma$, and that $\Delta'$ is a (finite,
geometric) \emph{simplicial subdivision} of $\Delta$, if (a) every simplex
of $\Sigma'$ is contained in some simplex of $\Sigma$; and (b) the union of
the simplices of $\Sigma'$ is equal to the union of the simplices of
$\Sigma$. Then, given any simplex $L \in \Sigma$ with corresponding face
$F \in \Delta$, the subcomplex $\Sigma'_L$ of $\Sigma'$ consisting of all
simplices of $\Sigma'$ contained in $L$ is called the \emph{restriction} of
$\Sigma'$ to $L$. The subcomplex $\Delta'_F$ of $\Delta'$ corresponding to
$\Sigma'_L$ is the \emph{restriction} of $\Delta'$ to $F$. Clearly,
$\Delta'_F$ is a simplicial subdivision of the abstract simplex $2^F$.

The more general notion of topological subdivision introduced in
\cite[Section~2]{Sta92} will not be needed here.

\bigskip
\noindent
\textbf{Barycentric and edgewise subdivisions.}
The (first, simplicial) \emph{barycentric subdivision} of a simplicial complex
$\Delta$, denoted $\sd (\Delta)$, is defined as the simplicial complex
consisting of all (possibly empty) chains $F_0 \subset F_1 \cdots \subset F_k$
of nonempty faces of $\Delta$. It is well known that $\sd (\Delta)$ is a
simplicial subdivision of $\Delta$. The restriction of $\sd (\Delta)$ to a
face $F \in \Delta$ coincides with the barycentric subdivision of $2^F$.

To define the $r$th edgewise subdivision of a simplicial complex, we follow
the discussions in \cite[Section~4]{BrW09} \cite[Section~6]{BR05}. Suppose that
$\Delta$ is a simplicial complex on the vertex set $\Omega_1 = \{ \ree_1,
\ree_2,\dots,\ree_m \}$ of coordinate vectors in some Euclidean space $\RR^m$.
For $\raa = (a_1, a_2,\dots,a_m) \in \ZZ^m$ we will denote
by $\supp(\raa)$ the set of indices $i \in \{1, 2,\dots,m\}$ for which $a_i
\ne 0$ and write $\iota(\raa) = (a_1, a_1 + a_2,\dots,a_1 + a_2 + \cdots +
a_m)$. The $r$th edgewise subdivision of $\Delta$ is the simplicial complex
$\Delta^{\langle r \rangle}$ on the ground set $\Omega_r = \{ (i_1,
i_2,\dots,i_m) \in \NN^m: i_1 + i_2 + \cdots + i_m = r\}$ of which a set $G \subseteq \Omega_r$ is a face if the following two conditions are satisfied:

\begin{itemize}
\itemsep=0pt
\item[$\bullet$]
$\bigcup_{\ruu \in G} \, \supp(\ruu) \in \Delta$ and

\item[$\bullet$]
$\iota(\ruu) - \iota(\rvv) \in \{0, 1\}^m$, or $\iota(\rvv) - \iota
(\ruu) \in \{0, 1\}^m$, for all $\ruu, \rvv \in G$.
\end{itemize}

\noindent
We leave it to the reader to verify that $\Delta^{\langle r \rangle} = \Delta$ for $r=1$ and that if $\Delta$ is a flag simplicial complex (meaning that every minimal non-face has at most two elements), then so is $\Delta^{\langle r \rangle}$.

For any $r \ge 1$, the simplicial complex $\Delta^{\langle r \rangle}$ can be realized as a simplicial subdivision of $\Delta$ (see, for instance, \cite[Section~6]{BR05}). The restriction of 
$\Delta^{\langle r \rangle}$ to a face $F \in \Delta$ coincides with the $r$th edgewise subdivision $(2^F)^{\langle r \rangle}$ of the simplex $2^F$; it has exactly $r^{\dim(F)}$ faces of the same dimension as $F$.

\bigskip
\noindent
\textbf{Face enumeration.}
Consider an $(n-1)$-dimensional simplicial complex $\Delta$ and let $f_i (\Delta)$
denote the number of $i$-dimensional faces of $\Delta$. The \emph{$h$-polynomial}
of $\Delta$ is defined as
  $$ h(\Delta, x) \ = \ \sum_{i=0}^n \, f_{i-1} (\Delta) \, x^i (1-x)^{n-i}. $$

\noindent
For the importance of $h$-polynomials, the reader is referred to
\cite[Chapter~II]{StaCCA}.

\begin{example} \label{ex:baryedge-h} \rm
For the barycentric subdivision $\Delta = \sd(2^V)$ of an $(n-1)$-dimensional
simplex $2^V$, it is well known that $h(\Delta, x)$ is equal to the Eulerian
polynomial $A_n (x)$. For example, for $n=3$ (see Figure~\ref{fig:KK33}) we
have $f_{-1} (\Delta) = 1$, $f_0 (\Delta) = 7$, $f_1 (\Delta) = 12$ and $f_2(\Delta) = 6$ and hence $h(\Delta, x) = 1 + 4x + x^2$.

An explicit formula can be given for the $h$-polynomial of the $r$th
edgewise subdivision of any $(n-1)$-dimensional simplicial complex $\Delta$.
Indeed, combining \cite[Corollary~6.8]{BR05} with \cite[Corollary~1.2]{BrW09},
one gets

  \begin{equation} \label{eq:hedge}
    h(\Delta^{\langle r \rangle}, x) \ = \ {\rm E}_r \left( (1 + x + x^2 +
    \cdots + x^{r-1})^n \, h (\Delta, x) \right),
  \end{equation}

\noindent
where ${\rm E}_r: \RR[x] \to \RR[x]$ is the linear operator defined by setting
${\rm E}_r (x^k) = x^{k/r}$, if $k$ is divisible by $r$, and ${\rm E}_r (x^k) =
0$ otherwise.
\end{example}

Given a simplicial subdivision $\Gamma$ of an $(n-1)$-dimensional simplex $2^V$,
the \emph{local $h$-polynomial}, denoted $\ell_V (\Gamma, x)$, of $\Gamma$ is
defined \cite[Definition~2.1]{Sta92} by Equation~(\ref{eq:deflocalh}), where
$\Gamma_F$ is the restriction of $\Gamma$ to the face $F \in 2^V$. Then $\ell_V
(\Gamma, x)$ has degree at most $n-1$ (unless $V = \varnothing$, in which case
$\ell_V (\Gamma, x) = 1$) and nonnegative and symmetric coefficients, so that
$x^n \, \ell_V (\Gamma, 1/x) = \ell_V (\Gamma, x)$. For examples and further
properties of local $h$-polynomials, see \cite[Part~I]{Sta92}
\cite[Section~II.10]{StaCCA} and \cite{Ath12, AS12, Sav13}.

\section{Proof of Theorem~\ref{thm:edgebary}}
\label{sec:proof}

This section deduces Theorem~\ref{thm:edgebary} from the results of
Section~\ref{sec:color} and Theorem~\ref{thm:dnr}.

\medskip
\noindent
\emph{Proof of Theorem~\ref{thm:edgebary}.} We set $\Gamma = \sd(2^V)^{\langle
r \rangle}$ and consider a $k$-element set $F \subseteq V$. As the relevant
discussions in Section~\ref{sec:sub} show, the restriction of $\Gamma$ to $F$
satisfies $\Gamma_F = \sd(2^F)^{\langle r \rangle}$ and $h(\sd(2^F), x) = A_k
(x)$. Thus, combining (\ref{eq:hedge}) with Proposition~\ref{prop:fdesfexc}
we get

  \begin{eqnarray*}
    h (\Gamma_F, x) &=& h (\sd(2^F)^{\langle r \rangle}, x) \ = \ {\rm E}_r
    \left( (1 + x + x^2 + \cdots + x^{r-1})^k \, A_k (x) \right) \\
    &=& {\rm E}_r \left( \sum_{w \in \ZZ_r \wr \mathfrak{S}_k} x^{\fex (w)}
    \right) \ = \ \sum_{w \in (\ZZ_r \wr \mathfrak{S}_k)^b} x^{\fex (w) / r},
  \end{eqnarray*}

\bigskip
\noindent
where $(\ZZ_r \wr \mathfrak{S}_k)^b$ is the set of balanced elements of $\ZZ_r
\wr \mathfrak{S}_k$. Since adding or removing fixed points from $w \in \ZZ_r
\wr \mathfrak{S}_n$ does not affect $\fex(w)$, the previous equalities, the
defining equation (\ref{eq:deflocalh}) and an easy application of the principle
of Inclusion-Exclusion show that

  $$ \ell_V (\Gamma, x) \ = \ \sum_{k=0}^n \, (-1)^{n-k} {n \choose k} \sum_{w
     \in (\ZZ_r \wr \mathfrak{S}_k)^b} x^{\fex (w) / r} \ = \
     \sum_{w \in (\dD^r_n)^b} x^{\fex (w) / r} $$

\bigskip
\noindent
and the proof of (\ref{eq:localedgebary}) follows.

We now assume that $r \ge 2$ and derive Equation~(\ref{eq:localgammaedgebary})
from Theorem~\ref{thm:dnr}. Since $d_n^r(x)$ has degree $n$ and zero constant
term, it can be written uniquely in the form
  \begin{equation} \label{eq:dnrsymsum}
    d_n^r (x) \ = \ f^+_{n,r} (x) \, + \, f^-_{n,r} (x),
  \end{equation}
where $f^+_{n,r} (x)$ and $f^-_{n,r} (x)$ are polynomials of degrees at most
$n-1$ and $n$, respectively, satisfying
  \begin{eqnarray}
  \label{eq:fnr+}  f^+_{n,r} (x) & = & x^n f^+_{n,r} (1/x) \\
  \label{eq:fnr-}  f^-_{n,r} (x) & = & x^{n+1} f^-_{n,r} (1/x)
\end{eqnarray}
(see, for instance, \cite[Lemma 2.4]{BeSt10} for this elementary fact). From
Proposition~\ref{prop:dnr1} we get that
  \begin{equation} \label{eq:dnrdecomp}
    d_n^r (x) \ = \ \sum_{w \in (\dD^r_n)^b} x^{\fex (w) / r} \ \ +
    \sum_{w \in \dD^r_n \sm (\dD^r_n)^b} x^{\lceil \frac{\fex(w)}{r} \rceil}.
  \end{equation}
It is an easy consequence of Proposition~\ref{prop:dnrflag} that the two
summands in the right-hand side of (\ref{eq:dnrdecomp}) have degrees at most
$n-1$ and $n$ and satisfy (\ref{eq:fnr+}) and (\ref{eq:fnr-}), respectively.
Thus, the uniqueness of the decomposition (\ref{eq:dnrsymsum}) implies that
  $$ f^+_{n,r} (x) \ = \ \sum_{w \in (\dD^r_n)^b} x^{\fex (w) / r}. $$
For similar reasons, the decomposition of $d_n^r(x)$ provided by
Theorem~\ref{thm:dnr} implies that
  $$ f^+_{n,r} (x) \ = \ \sum_{i=0}^{\lfloor n/2 \rfloor} \xi^+_{n, r, i} \,
     x^i (1+x)^{n-2i} $$
and the proof of (\ref{eq:localgammaedgebary}) follows.
\qed

\section{Remarks and open problems}
\label{sec:rem}

\noindent
1. Chow and Mansour showed \cite[Theorem~5~(ii)]{CM10} that $d^r_n (x)$ is
real-rooted and, as a result, that it has log-concave and unimodal coefficients.
Theorem~\ref{thm:dnr} provides a transparent proof of the unimodality of $d^r_n
(x)$. It also shows that its peak occurs at $\lfloor (n+1)/2 \rfloor$.

\bigskip
\noindent
2. Following Steingr\'imsson~\cite[Definition~36]{Stei94}, we call a colored
permutation $w \in \ZZ_r \wr \mathfrak{S}_n$ \emph{reverse alternating} if for
$i \in \{1, 2,\dots,n\}$ we have: $i$ is a descent of $w$ if and only if $i$
is odd. For fixed $r$, the exponential generating function of these permutations
was computed in \cite[Theorem~39]{Stei94}. The following statement specializes
to \cite[Theorem~11]{Zh13} for $r=2$.

\begin{corollary} \label{cor:dnr}
For all positive integers $n, r$, the number $(-1)^{\lfloor \frac{n+1}{2}
\rfloor} \, d^r_n (-1)$ is equal to the number of reverse alternating
colored permutations in $\ZZ_r \wr \mathfrak{S}_n$.
\end{corollary}

\noindent
\emph{Proof.} Setting $x = -1$ in (\ref{eq:dnrdecompose}) we get

  $$ d^r_n (-1) \ = \ \cases{ (-1)^{\frac{n}{2}} \ \xi^+_{n, r, \frac{n}{2}},
                              & if \ $n$ is even, \cr
                              (-1)^{\frac{n+1}{2}} \ \xi^-_{n, r, \frac{n+1}{2}},
                              & if \ $n$ is odd} $$

\noindent
and the proof follows from the combinatorial interpretations given to the
$\xi^+_{n, r, i}$ and $\xi^-_{n, r, i}$ in Theorems~\ref{thm:bary} and
\ref{thm:edgebary}.
\qed

\bigskip
\noindent
3. The proof of Theorem~\ref{thm:edgebary} shows that the polynomials
$d_n^2 (x)$ and $f_n^2 (x)$ can be computed one from the other in a simple
way. Indeed, (\ref{eq:dnrdecomp}) shows that $d_n^2 (x) = \alpha_n (x) +
\beta_n (x)$, where $\alpha_n (x)$ and $\beta_n (x)$ satisfy $f_n^2 (x) =
\alpha_n (x^2) + \beta_n (x^2) / x$ and hence are uniquely determined by
$f_n^2 (x)$. Similarly, $\alpha_n (x)$ and $\beta_n (x)$ are uniquely
determined by $d_n^2 (x)$ by the uniqueness of the decomposition
(\ref{eq:dnrsymsum}), already discussed.

\bigskip
\noindent
4. Does the local $h$-polynomial of the subdivision $\sd(2^V)^{\langle r \rangle}$, computed in Theorem~\ref{thm:edgebary}, have (when not identically zero) only real roots?

Since we are not aware of an explicit example of a flag (geometric) simplicial subdivision of the simplex whose local $h$-polynomial is not real-rooted, even a negative answer to this question would be of interest. The answer is affirmative for $r=1$ by the result of 
\cite{Zh95} and, as was verified by Savvidou~\cite{Sav13b} by explicit
computation, for $r = 2$ and $n \le 50$.


\begin{thebibliography}{99}
%
\bibitem{ABR01}
R.M.~Adin, F.~Brenti and Y.~Roichman,
\textit{Descent numbers and major indices for the hyperoctahedral group},
Adv. in Appl. Math. {\bf~27} (2001), 210--224.
%
\bibitem{Ath12}
C.A.~Athanasiadis,
\textit{Flag subdivisions and $\gamma$-vectors},
Pacific J. Math. {\bf~259} (2012), 257--278.
%
\bibitem{AS12}
C.A.~Athanasiadis and C.~Savvidou,
\textit{The local $h$-vector of the cluster subdivision of a simplex},
S\'em. Lothar. Combin. {\bf~66} (2012), Article B66c, 21pp (electronic).
%
\bibitem{BB07}
E.~Bagno and R.~Biagioli,
\textit{Colored descent representations for complex reflection groups},
Israel J. Math. {\bf~160} (2007), 317--347.
%
\bibitem{BG06}
E.~Bagno and D.~Garber,
\textit{On the excedance number of colored permutation groups},
S\'em. Lothar. Combin. {\bf~53} (2006), Article B53f, 17pp (electronic).
%
\bibitem{BeSt10}
M.~Beck and A.~Stapledon,
\textit{On the log-concavity of Hilbert series of Veronese subrings and
Ehrhart series},
Math. Z. {\bf~264} (2010), 195--207.
%
\bibitem{Bj95}
A.~Bj\"orner,
\textit{Topological methods}, in \textit{Handbook of combinatorics}
(R.L.~Graham, M.~Gr\"otschel and L.~Lov\'asz, eds.),
North Holland, Amsterdam, 1995, pp.~1819--1872.
%
\bibitem{BjW05}
A.~Bj\"orner and V.~Welker,
\textit{Segre and Rees products of posets, with ring-theoretic applications},
J. Pure Appl. Algebra {\bf~198} (2005), 43--55.
%
\bibitem{Bre90}
F.~Brenti,
\textit{Unimodal polynomials arising from symmetric functions},
Proc. Amer. Math. Soc. {\bf~108} (1990), 1133--1141.
%
\bibitem{BrW09}
F.~Brenti and V.~Welker,
\textit{The Veronese construction for formal power series and graded algebras},
Adv. in Appl. Math. {\bf~42} (2009), 545--556.
%
\bibitem{BR05}
M.~Brun and T.~R\"omer,
\textit{Subdivisions of toric complexes},
J. Algebraic Combin. {\bf~21} (2005), 423--448.
%
\bibitem{CTZ09}
W.Y.C.~Chen, R.L.~Tang and A.F.Y.~Zhao,
\textit{Derangement polynomials and excedances of type $B$},
Electron. J. Combin. {\bf~16} (2) (2009), Research Paper 15, 16pp (electronic).
%
\bibitem{Ch09}
C.-O.~Chow,
\textit{On derangement polynomials of type $B$. II},
J. Combin. Theory Series A {\bf~116} (2009), 816--830.
%
\bibitem{CM10}
C.-O.~Chow and T.~Mansour,
\textit{Counting derangements, involutions and unimodal elements in the wreath
product $C_r \wr \mathfrak{S}_n$},
Israel J. Math. {\bf~179} (2010), 425--448.
%
\bibitem{EG00}
H.~Edelsbrunner and D.R.~Grayson,
\textit{Edgewise subdivision of a simplex},
Discrete Comput. Geom. {\bf~24} (2000), 707--719.
%
\bibitem{FH09}
D.~Foata and G.-N.~Han,
\textit{Signed words and permutations, V; a sextuple distribution},
Ramanujan J. {\bf~19} (2009), 29--52.
%
\bibitem{FH11}
D.~Foata and G.-N.~Han,
\textit{The decrease value theorem with an application to permutation statistics},
Adv. in Appl. Math. {\bf~46} (2011), 296--311.
%
\bibitem{Gr89}
D.R.~Grayson,
\textit{Exterior power operations on higher $K$-theory},
K Theory {\bf~3} (1989), 247--260.
%
\bibitem{LSW12}
S.~Linusson, J.~Shareshian and M.L.~Wachs,
\textit{Rees products and lexicographic shellability},
J. Combinatorics {\bf~3} (2012), 243--276.
%
\bibitem{Mat03}
J.~Matousek,
Using the Borsuk-Ulam Theorem: Lectures on Topological Methods in
Combinatorics and Geometry,
Universitext, Springer, 2003.
%
\bibitem{Mo13}
P.~Mongelli,
\textit{Excedances in classical and affine Weyl groups},
J. Combin. Theory Series A {\bf~120} (2013), 1216--1234.
%
\bibitem{Sav13}
C.~Savvidou,
\textit{Barycentric subdivisions, clusters and permutation enumeration} (in Greek),
Doctoral Dissertation, University of Athens, 2013.
%
\bibitem{Sav13b}
C.~Savvidou,
\textit{personal communication}, 2013.
%
\bibitem{SW09}
J.~Shareshian and M.L.~Wachs,
\textit{Poset homology of Rees products and $q$-Eulerian polynomials},
Electron. J. Combin. {\bf~16} (2) (2009), Research Paper 20, 29pp (electronic).
%
\bibitem{SW10}
J.~Shareshian and M.L.~Wachs,
\textit{Eulerian quasisymmetric functions},
Adv. Math. {\bf~225} (2010), 2921--2966.
%
\bibitem{SZ12}
H.~Shin and J.~Zeng,
\textit{The symmetric and unimodal expansion of Eulerian polynomials via continued 
fractions},
European J. Combin. {\bf~33} (2012), 111--127.
%
\bibitem{Sta92}
R.P.~Stanley,
\textit{Subdivisions and local $h$-vectors},
J. Amer. Math. Soc. {\bf~5} (1992), 805--851.
%
\bibitem{StaCCA}
R.P.~Stanley,
Combinatorics and Commutative Algebra,
second edition, Birkh\"auser, Basel, 1996.
%
\bibitem{StaEC1}
R.P.~Stanley,
Enumerative Combinatorics, vol.~1,
Cambridge Studies in Advanced Mathematics {\bf~49},
Cambridge University Press, second edition, Cambridge, 2011.
%
\bibitem{Stei92}
E.~Steingr\'imsson,
\textit{Permutation statistics of indexed and poset permutations},
Ph.D. thesis, MIT, 1992.
%
\bibitem{Stei94}
E.~Steingr\'imsson,
\textit{Permutation statistics of indexed permutations},
European J. Combin. {\bf~15} (1994), 187--205.
%
\bibitem{Wa07}
M.~Wachs,
\textit{Poset Topology: Tools and Applications}, in
\emph{Geometric Combinatorics} (E.~Miller, V.~Reiner and B.~Sturmfels, eds.),
IAS/Park City Mathematics Series {\bf~13},
pp.~497--615, Amer. Math. Society, Providence, RI, 2007.
%
\bibitem{Zh95}
X.~Zhang,
\textit{On $q$-derangement polynomials},
in \textit{Combinatorics and Graph Theory '95}, Vol. 1 (Hefei),
pp.~462--465, World Sci. Publishing, River Edge, NJ, 1995.
%
\bibitem{Zh13}
A.F.Y.~Zhao,
\textit{Excedance numbers for the permutations of type $B$},
Electron. J. Combin. {\bf~20} (2) (2013), Research Paper 28, 20pp (electronic).
%
\end{thebibliography}
\end{document}